\documentclass[reqno,a4paper,10pt]{article}
\usepackage{color}
\usepackage[latin1]{inputenc}
\usepackage[english]{babel}
\usepackage{amsmath,amssymb,amsfonts,amsthm,enumerate}
\usepackage{mathrsfs}
\usepackage[T1]{fontenc}
\usepackage{epsfig}
\usepackage{graphicx}
\usepackage[
color,
final] 
{showkeys}

\definecolor{refkeybis}{gray}{.65}
\definecolor{labelkeybis}{gray}{.65}
{\makeatletter
\def\SK@refcolor{\color{refkeybis}}%
\def\SK@labelcolor{\color{labelkeybis}}}

\newtheorem{theorem}{Theorem}[section]

\newtheorem{definition}[theorem]{Definition}
\newtheorem{remark}[theorem]{Remark}

\newtheorem{example}[theorem]{Example}

\newcommand{\R}{\mathbb{R}}
\newcommand{\C}{\mathbb{C}}

\newcommand{\Leb}[1]{{\mathscr L}^{#1}} 

\renewcommand{\P}{\mathbb{P}}

\newcommand{\e}{\varepsilon}

\newcommand{\n}{\nabla}

\newcommand{\p}{\partial}

\newcommand{\Probabilities}[1]{\mathscr P\bigl(#1\bigr)} 
\newcommand{\Measuresp}[1]{\mathscr M_+\bigl(#1\bigr)} 

\newcommand{\E}{{\mathbb E}}

 \newcommand{\bb}{{\mbox{\boldmath$b$}}}
 \newcommand{\cc}{{\mbox{\boldmath$c$}}}

 \newcommand{\tauV}{{\kern-3pt\tau}}

 \newcommand{\sxX}{{\mbox{\scriptsize\boldmath$X$}}}

 \newcommand{\XX}{{\mbox{\boldmath$X$}}}
 \newcommand{\YY}{{\mbox{\boldmath$Y$}}}
 
 \newcommand{\oVVVk}{\overline{\mbox{\boldmath$V$}}\kern-3pt}
 \newcommand{\tVVVk}{\tilde{\mbox{\boldmath$V$}}\kern-3pt}

 \newcommand{\mmu}{{\mbox{\boldmath$\mu$}}}
 \newcommand{\nnu}{{\mbox{\boldmath$\nu$}}}

 \newcommand{\eps}{\varepsilon}

\title{Almost everywhere well-posedness of continuity equations\\
with measure initial data}
\date{}
\author{Luigi Ambrosio\
   \thanks{\textsf{l.ambrosio@sns.it}}
   \and
   Alessio Figalli\
   \thanks{\textsf{figalli@math.utexas.edu}}}

\begin{document}

\maketitle
\begin{center}
{\bf Abstract}
\end{center}
{\small The aim of this note is to present some new results concerning ``almost everywhere''
well-posedness and stability of continuity equations with measure initial data.
The proofs of all such results can be found in \cite{amfifrgi}, together
with some application to the semiclassical limit of the Schr\"odinger equation.
}
\begin{center}
{\bf Resum\'e}
\end{center}
{\small Dans cette note, nous pr\'esentons des nouveaux r\'esultats concernant
l'existence, l'unicit\'e (au sens ``presque partout'') et la stabilit\'e pour des \'equations
de continuit\'e avec donn\'ees initiales mesures.
Les preuves de tous ces r\'esultats sont donn\'ees dans \cite{amfifrgi}, avec aussi des applications
\`a la limite semiclassique pour l'\'equation de Schr\"odinger.
}
\bigskip

Starting from the seminal paper of DiPerna-Lions \cite{diplions}
(dealing mostly with the transport equation), in
\cite{ambrosio,cetraro} the well-posedness of the continuity
equation
\begin{equation}\label{contieq}
\frac{d}{dt}\mu_t+\nabla\cdot (\bb_t\mu_t)=0.
\end{equation}
has been stronly related to well-posedness of the ODE (here we use
the notation $\bb(t,x)=\bb_t(x)$)
\begin{equation}\label{ODE}
\left\{\begin{array}{ll}
\dot \XX(t,x)=\bb_t(\XX(t,x))&\text{for $\Leb{1}$-a.e. $t\in (0,T)$,}\\
\XX(0,x)=x,
\end{array}
\right.
\end{equation}
for ``almost every'' $x\in\R^d$.
More precisely, observe the concept of solution to \eqref{ODE} is not
invariant under modification of $\bb$ in Lebesgue negligible
sets, while many applications of the theory to fluid dynamics (see
for instance \cite{lions2}, \cite{lions3}) and conservation laws
need this invariance property. This leads to the concept of
\emph{regular Lagrangian flow} (RLF in short): one may ask that,
for all $t\in [0,T]$, the image $\XX(t,\cdot)_\sharp\Leb{d}$
of the Lebesgue measure $\Leb{d}$ under the
flow map $x\mapsto\XX(t,x)$ is still controlled by $\Leb{d}$ (see
Definition~\ref{RLflow} below). Then, existence and uniqueness (up to $\Leb{d}$-negligible sets)
and stability of the RLF $\XX(t,x)$ in $\R^d$ hold true provided
the functional version of \eqref{contieq}, namely
\begin{equation}\label{contieqw}
\frac{d}{dt}w_t+\nabla\cdot (\bb_t w_t)=0,
\end{equation}
is well-posed in the set of non-negative bounded integrable funtions
$L^\infty_+\bigl([0,T];L^1(\R^d)\cap L^\infty(\R^d)\bigr)$.

Now, we may view \eqref{contieq} as an infinite-dimensional ODE in
$\Probabilities{\R^d}$, the space of probability measures in
$\R^d$ and try to obtain existence and uniqueness results for \eqref{contieq}
in the same spirit of the finite-dimensional theory, starting from the
simple observation that $t\mapsto\delta_{\sxX(t,x)}$ solves \eqref{contieq}.
We may expect that
if we fix a ``good'' measure $\nnu$ in the space $\Probabilities{\R^d}$
of initial data, then existence, uniqueness $\nnu$-a.e. and stability hold.
Moreover, for $\nnu$-a.e. $\mu$, the unique and stable
solution of \eqref{contieq} starting from $\mu$
should be given by
\begin{equation}\label{ovvia}
\mmu(t,\mu):=\int \delta_{\sxX(t,x)}\,d\mu(x)\qquad \forall\, t\in
[0,T],\,\,\mu\in\Probabilities{\R^d}.
\end{equation}

\section{Continuity equations and flows}

We use a standard and hopefully self-explainatory notation.
Let $\bb:[0,T]\times\R^d\to\R^d$ be a Borel vector field belonging to
$L^1_{\rm loc}\bigl([0,T]\times\R^d\bigr)$, and
set $\bb_t(\cdot):=\bb(t,\cdot)$; we \emph{shall not} work with the
Lebesgue equivalence class of $\bb$, although a posteriori the
theory is independent of the choice of the representative.

\begin{definition}[$\nu$-RLF in $\R^d$]\label{RLflow}
Let $\XX(t,x):[0,T]\times\R^d\to\R^d$ and
$\nu\in {\mathscr M}_+(\R^d)$ with $\nu\ll\Leb{d}$ and with bounded
density. We say that $\XX(t,x)$ is a $\nu$-RLF in $\R^d$ (relative
to $\bb$) if the following two conditions are fulfilled:
\begin{itemize}
\item[(i)] for $\nu$-a.e. $x$, the function $t\mapsto\XX(t,x)$
is an absolutely continuous integral solution to the ODE \eqref{ODE}
in $[0,T]$ with $\XX(0,x)=x$;
\item[(ii)] $\XX(t,\cdot)_\sharp\nu\leq C\Leb{d}$ for all $t\in
[0,T]$, for some constant $C$ independent of $t$.
\end{itemize}
\end{definition}

By a simple application of Fubini's theorem this concept is, unlike the
single condition (i), invariant in the Lebesgue equivalence class of $\bb$.
In this context, since all admissible initial measures $\nu$ are
bounded above by $C\Leb{d}$, uniqueness of the $\nu$-RLF can and
will be understood in the following stronger sense: if $f,\,g\in
L^1(\R^d)\cap L^\infty(\R^d)$ are nonnegative and $\XX$ and $\YY$
are respectively a $f\Leb{d}$-RLF and a $g\Leb{d}$-RLF, then
$\XX(\cdot,x)=\YY(\cdot,x)$ for $\Leb{d}$-a.e. $x\in\{f>0\}\cap
\{g>0\}$.

\begin{remark}{\rm We recall that the $\nu$-RLF exists
for all $\nu\leq C\Leb{d}$, and is unique, in the strong sense
described above, under the following
assumptions on $\bb$: $|\bb|$ is uniformly bounded, $\bb_t\in
BV_{\rm loc}(\R^d;\R^d)$ and $\nabla\cdot\bb_t=g_t\Leb{d}\ll\Leb{d}$
for $\Leb{1}$-a.e. $t\in (0,T)$, with
$$
\|g_t\|_{L^\infty(\R^d)}\in L^1(0,T),\qquad |D\bb_t|(B_R)\in
L^1(0,T)\quad\text{for all $R>0$,}
$$
where $|D\bb_t|$ denotes the total variation of the distributional
derivative of $\bb_t$. (See \cite{ambrosio} or
\cite{cetraro} and the paper \cite{bouchut} for Hamiltonian vector
fields.)}\end{remark}

Given a nonnegative $\sigma$-finite measure
$\nnu\in\Measuresp{\Probabilities{\R^d}}$, we denote by
$\E\nnu\in\Measuresp{\R^d}$ its expectation, namely
$$
\int_{\R^d}\phi\,d\E\nnu=\int_{{\mathscr
P}(\R^d)}\int_{\R^d}\phi\,d\mu \,d\nnu(\mu)\qquad\text{for all
$\phi$ bounded Borel.}
$$

\begin{definition}[Regular measures in
$\Measuresp{\Probabilities{\R^d}}$]\label{RegMis} Let
$\nnu\in\Measuresp{\Probabilities{\R^d}}$. We say that $\nnu$ is
\emph{regular} if $\E\nnu\leq C\Leb{d}$ for some constant $C$.
\end{definition}

\begin{example}\label{eRegMis}
{\rm (1) The first standard example of a regular measure $\nnu$ is
the law under $\rho\Leb{d}$ of the map $x\mapsto\delta_x$, with
$\rho\in L^1(\R^{d})\cap L^\infty(\R^d)$ nonnegative. Actually, one can
even consider the law under $\Leb{d}$, and in this case
$\nnu$ would be $\sigma$-finite instead of finite.

(2) If $d=2n$ and $z=(x,p)\in\R^n\times\R^n$ (this factorization
corresponds for instance to flows in a phase space), instead of
considering the law of under $\rho\Leb{2n}$ of the map
$(x,p)\mapsto\delta_x\otimes\delta_p$, one may also consider the law
under $\rho\Leb{n}$ of the map $x\mapsto \delta_x\times\gamma$, with
$\rho\in L^1(\R^n_x)\cap L^\infty(\R^n_x)$ nonnegative and
$\gamma\in\Probabilities{\R^n_p}$ bounded from above by a constant
multiple of $\Leb{n}$.}
\end{example}

We observe that Definition~\ref{RLflow} has a
natural (but not perfect) transposition to flows in
$\Probabilities{\R^d}$:

\begin{definition}[Regular Lagrangian flow in
$\Probabilities{\R^d}$]\label{RLflowmis} Let
$\mmu:[0,T]\times\Probabilities{\R^d}\to\Probabilities{\R^d}$ and
$\nnu\in\Measuresp{\Probabilities{\R^d}}$. We say that $\mmu$ is a
$\nnu$-RLF in $\Probabilities{\R^d}$ (relative to $\bb$) if
\begin{itemize}
\item[(i)] for $\nnu$-a.e. $\mu$, $|\bb|\in L^1_{\rm
loc}\bigl((0,T)\times\R^d;\mu_tdt\bigr)$,
$t\mapsto\mu_t:=\mmu(t,\mu)$ is continuous from $[0,1]$ to
$\Probabilities{\R^d}$ with $\mmu(0,\mu)=\mu$ and $\mu_t$ solves
\eqref{contieq} in the sense of distributions;
\item[(ii)] $\E(\mmu(t,\cdot)_\sharp\nnu)\leq C\Leb{d}$ for all $t\in [0,T]$, for
some constant $C$ independent of $t$.
\end{itemize}
\end{definition}

Notice that condition (ii) is weaker than $\mmu(t,\cdot)_\sharp\nnu\leq C\nnu$
(which would be the analogue of (ii) in Definition~\ref{RLflow} if
we were allowed to choose $\nu=\Leb{d}$), and it is actually sufficient
and much more flexible for our purposes, since we would like to consider
measures $\nnu$ generated as in Example~\ref{eRegMis}(2).

\section{Existence, uniqueness and stability of the RLF}

In this section we recall the main existence and uniqueness
results of the $\nu$-RLF in $\R^d$, and see their extensions to
$\nnu$-RLF in $\Probabilities{\R^d}$. The following result is proved
in \cite[Theorem~19]{cetraro} for the
part concerning existence and in \cite[Theorem~16,
Remark~17]{cetraro} for the part concerning uniqueness.

\begin{theorem}[Existence and uniqueness of the $\nu$-RLF in
$\R^d$]\label{texirlfrd} Assume that \eqref{contieqw} has existence
and uniqueness in $L^\infty_+\bigl([0,T];L^1(\R^d)\cap
L^\infty(\R^d)\bigr)$. Then, for all $\nu\ll\Leb{d}$ with bounded
density the $\nu$-RLF exists and is unique.
\end{theorem}

The next result shows that, uniqueness of \eqref{contieqw} in
$L^\infty_+\bigl([0,T];L^1(\R^d)\cap L^\infty(\R^d)\bigr)$ implies a
stronger property, namely uniqueness of the $\nnu$-RLF.

\begin{theorem}[Existence and uniqueness of the $\nnu$-RLF in
$\Probabilities{\R^d}$]\label{texirlfprob} Assume that
\eqref{contieqw} has uniqueness in $L^\infty_+\bigl([0,T];L^1(\R^d)\cap
L^\infty(\R^d)\bigr)$. Then, for all
$\nnu\in\Measuresp{\Probabilities{\R^d}}$ regular,
there exists at most one $\nnu$-RLF in $\Probabilities{\R^d}$.
If \eqref{contieqw} has existence in  $L^\infty_+\bigl([0,T];L^1(\R^d)\cap
L^\infty(\R^d)\bigr)$, this unique flow is given by
\begin{equation}\label{realRLF}
\mmu(t,\mu):=\int_{\R^d}\delta_{\sxX(t,x)}\,d\mu(x),
\end{equation}
where $\XX(t,x)$ denotes the unique $\E\nnu$-RLF.
\end{theorem}

For the applications it is important to show that RLF's not only
exist and are unique, but also that they are stable. In the
statement of the stability result we shall consider measures
$\nnu_n\in\Probabilities{\Probabilities{\R^d}}$, $n\geq 1$, and a
limit measure $\nnu$. We shall assume that $\nnu_n=(i_n)_\sharp\P$,
where $(W,{\mathcal F},\P)$ is a probability measure space and
$i_n:W\to\Probabilities{\R^d}$ are measurable; we shall also assume
that $\nnu=i_\sharp\P$, with $i_n\to i$ $\P$-almost everywhere.
(Recall that Skorokhod theorem (see \cite[\S8.5, Vol.
II]{bogachevII}) shows that weak convergence of $\nnu_n$ to $\nnu$
always implies this sort of representation, even with $W=[0,1]$
endowed with the standard measure structure, for suitable
$i_n,\,i$.) The following formulation of the stability result is
particularly suitable for the application to semiclassical limit of
the Schr\"odinger equation.

Henceforth, we fix an autonomous vector field $\bb:\R^d\to\R^d$
satisfying the following regularity conditions:

\begin{itemize}
\item[(a)] $d=2n$ and $\bb(x,p)=(p,\cc(x))$, $(x,p)\in\R^d$,
$\cc:\R^n\to\R^n$ Borel and locally integrable;
\item[(b)] there
exists a closed $\Leb{n}$-negligible set $S$ such that $\cc$ is
locally bounded on $\R^n\setminus S$.
\end{itemize}

\begin{theorem}[Stability of the $\nnu$-RLF in
$\Probabilities{\R^d}$]\label{tstable} Let $i_n,\,i$ be as above and
let $\mmu_n:[0,T]\times i_n(W)\to\Probabilities{\R^d}$ be satisfying
$\mmu_n(0,i_n(w))=i_n(w)$ and the following conditions:
\begin{itemize}
\item[(i)] (uniform regularity)
$$\sup_{n\geq 1}\sup_{t\in [0,T]}\int_W\int_{\R^d}\phi\,d\mmu_n(t,i_n(w))\,d\P(w)\leq
C\int_{\R^d} \phi\,dx
$$
for all $\phi\in C_c(\R^d)$ nonnegative;
\item[(ii)] (uniform decay away from $S$) for some $\beta>1$
\begin{equation}\label{ali4}
\sup_{\delta>0}\limsup_{n\to\infty}\int_W\int_0^T
\int_{B_R}\frac{1}{{\rm
dist}^\beta(x,S)+\delta}\,d\mmu_n(t,i_n(w))\,dt \,d\P(w)<\infty
\qquad\forall\, R>0;
\end{equation}
\item[(iii)] (space tightness) for all $\eps>0$,
$\P\Bigl(\bigl\{w:\ \sup\limits_{t\in
[0,T]}\mmu_n(t,i_n(w))(\R^d\setminus B_R)>\eps\bigr\}\Bigr)\to 0$
as $R\to\infty$;
\item[(iv)] (time tightness) for $\P$-a.e. $w\in
W$, for all $n\geq 1$ and $\phi\in C^\infty_c(\R^d)$,
$t\mapsto\int_{\R^d}\phi\,d\mmu_n(t,i_n(w))$ is absolutely
continuous in $[0,T]$ and
$$
\lim_{M\uparrow\infty}\P\biggl(\Bigl\{w\in W:\
\int_0^T\biggl|\biggr(\int_{\R^d}\phi\,d\mmu_n(t,i_n(w))\biggr)'\biggr|\,dt>M\Bigr\}\biggr)=0;
$$
\item[(v)] (limit continuity equation)
\begin{equation}\label{ali6}
\lim_{n\to\infty}\int_W
\biggl|\int_0^T\biggl[\varphi'(t)\int_{\R^d}\phi\,d\mmu_n(t,i_n(w))+\varphi(t)\int_{\R^d}
\langle\bb,\nabla\phi\rangle\,d\mmu_n(t,i_n(w))\biggr]\,dt\biggr|\,d\P(w)=0
\end{equation}
for all $\phi\in C^\infty_c\bigl(\R^d\setminus(S\times\R^n)\bigr)$,
$\varphi\in C^\infty_c(0,T)$.
\end{itemize}
Assume, besides (a), (b) above, that \eqref{contieqw} has
uniqueness in $L^\infty_+\bigl([0,T];L^1\cap L^\infty(\R^d)\bigr)$.
Then the $\nnu$-RLF $\mmu(t,\mu)$ relative to $\bb$ exists, is
unique (by Theorem~\ref{texirlfprob}) and
\begin{equation}\label{cetraro1}
\lim_{n\to\infty} \int_W\sup_{t\in [0,T]}d_{{\mathscr
P}}(\mmu_n(t,i_n(w)),\mmu(t,i(w)))\,d\P(w)=0
\end{equation}
where $d_{{\mathscr P}}$ is any bounded distance in $\Probabilities{\R^d}$ inducing
weak convergence of measures.
\end{theorem}

An example of application of the above stability result is the following:
let $\alpha\in (0,1)$ and let $\psi^\e_{x_0,p_0}:[0,T]\times \R^n \to \C$
be a family of solutions to the Schr\"odinger equation
\begin{equation}\label{wkb}
\left\{
\begin{array}{l} i\e \p_t \psi_{x_0,p_0}^\e(t)=-\frac{\e^2}{2}\Delta \psi_{x_0,p_0}^\e(t)+U\psi_{x_0,p_0}^\e(t)\\
\psi^\e_{x_0,p_0}(0)=\e^{-n\alpha/2}\phi_0\Bigl(\frac{x-x_0}{\e^{\alpha}}\Bigr)e^{i(x\cdot
p_0)/\e},
\end{array}
\right.
\end{equation}
with $\phi_0\in C^2_c(\R^n)$ and $\int|\phi_0|^2\,dx=1$.
When the potential $U$  is of class $C^2$, it was proven in \cite{gerard, lionspaul}
that for every $(x_0,p_0)$ the Wigner transforms
$W_\e\psi^\e_{x_0,p_0}(t)$ converge, in the natural
dual space ${\cal A}'$ \cite{lionspaul} for the Wigner transforms,
to $\delta_{\sxX(t,x_0,p_0)}$ as $\e\downarrow  0$.
Here $\XX(t,x,p)$ if the unique flow in $\R^{2n}$
associated to the Liouville equation
\begin{equation}
\label{eq:liouville}
\p_t W + p\cdot \n_x W-\n U(x) \cdot \n_p W=0.
\end{equation}
In \cite{amfifrgi}, relying also on some a-priori estimates of
\cite{amfrgi} (see also \cite{figpaul}), the authors consider a
potential $U$ which can be written as the sum of a repulsive Coulomb
potential $U_s$ plus a bounded Lipschitz interation term $U_b$ with
$\n U_b\in BV_{\rm loc}$. We observe that in this case the equation
\eqref{eq:liouville} does not even make sense for measure initial
data, as $\n U$ is not continuous. Still, they can prove \emph{full}
convergence as $\e\downarrow 0$, namely
\begin{equation}\label{limeps}
\lim_{\e\downarrow 0}
\int_{\R^d}\rho(x_0,p_0)\sup_{t\in [-T,T]}
d_{{\cal A}'}\bigl(W_\e\psi^\e_{x_0,p_0}(t),\delta_{\sxX(t,x_0,p_0)}\bigr)
dx_0dp_0=0
\qquad\forall\, T>0
\end{equation}
for all $\rho\in L^1(\R^{2n})\cap L^\infty(\R^{2n})$ nonnegative,
where $\XX(t,x,p)$ if the unique $\Leb{2n}$-RLF
associated to \eqref{eq:liouville} and $d_{{\cal A}'}$ is a bounded distance
inducing the weak$^*$ topology in the unit ball of ${\cal A'}$.

The proof of \eqref{limeps} relies on an application of
Theorem~\ref{tstable} to the Husimi transforms of
$\psi^\e_{x_0,p_0}(t)$. The scheme is sufficiently flexible to allow
more general families of initial conditions displaying partial
concentration, of position or momentum, or no concentration at all:
for instance, the limiting case $\alpha=1$ in \eqref{wkb} (related to
Example~\ref{eRegMis}(2)) leads to
$$
\lim_{\e\downarrow 0} \int_{\R^d}\rho(x_0)\sup_{t\in [-T,T]}
d_{{\cal
A}'}\bigl(W_\e\psi^\e_{x_0,p_0}(t),\mmu(t,\mu(x_0,p_0))\bigr) dx_0=0
\qquad\forall\, p_0\in\R^n,\,T>0
$$
for all $\rho\in L^1(\R^n)\cap L^\infty(\R^n)$ nonnegative, with
$\mmu(t,\mu)$ given by \eqref{ovvia} and
$\mu(x_0,p_0)=\delta_{x_0}\times|\hat\phi_0|^2(\cdot-p_0)\Leb{n}$.

\end{document}